
\input amstex
\documentstyle{amsppt}

\magnification=\magstep1
\hsize=6.5truein
\vsize=9truein
\baselineskip5pt

 at 10truept
 at 10truept
\font \smallrm=cmr10 at 10truept
 at 7truept
\font \smallbf=cmbx10 at 10truept
 at 10truept
\font \smallsl=cmsl10 at 10truept


\def \smallcirc {\, {\scriptstyle \circ} \,}
\def \N {\Bbb N}
\def \Z {\Bbb Z}
\def \C {\Bbb C}
\def \Cq {{\Bbb C}(q)}

\def \rcalzero {{{\Cal R}^{(0)}}}

\def \runo {{R^{(1)}}}

\def \gerg {{\frak g}}
\def \gerh {{\frak h}}
\def \ghat {{\hat \gerg}}
\def \hhat {{\hat \gerh}}

\def \edot {\dot E}
\def \fdot {\dot F}
\def \echeck {\check E}
\def \fcheck {\check F}
\def \edotcheck {\,\dot {\!{\echeck}}}
\def \fdotcheck {\,\dot {\!{\fcheck}}}

\def \phipre {\Phi_+^{\text{re}}}
\def \phipim {\Phi_+^{\text{im}}}
\def \phitildep {\widetilde{\Phi}_+}
\def \phitildepim {\widetilde{\Phi}_+^{\text{im}}}

\def \uqg {U_q(\ghat)}
\def \calUqg {{\Cal U}_q(\ghat)}
\def \uqm {U_q^-}
\def \calUqm {{\Cal U}_q^-}
\def \uqz {U_q^0}

\def \uqp {U_q^+}

\def \uqbm {U_q^\leq}
\def \calUqbm {{\Cal U}_q^\leq}
\def \uqbp {U_q^\geq}
\def \calUqbp {{\Cal U}_q^\geq}


\document

\topmatter

{\ }

\vskip-37pt

 \centerline{ \smallrm {\smallsl Journal of Pure and Applied Algebra\/}  {\smallbf 155}  no.~1 (2001), 41--52 }
 \vskip3pt
 \centerline{ \smallrm DOI: 10.1016/S0022-4049(99)00117-6 }

\vskip33pt

\title
    The  $ R $--matrix  action of untwisted affine
quantum groups at roots of 1
\endtitle

\author
   Fabio Gavarini
\endauthor

  \rightheadtext{ $ R $--matrix  action of affine quantum
groups at roots of 1 }

\affil
   Universit\`a degli Studi di Roma ``Tor Vergata'' --- Dipartimento
di Matematica  \\
   Via della Ricerca Scientifica, I-00133 Roma --- ITALY  \\
\endaffil

\address\hskip-\parindent
           Universit\`a degli Studi di Roma ``Tor Vergata''  \newline
           Dipartimento di Matematica  \newline
           Via della Ricerca Scientifica  \newline
           I-00133 Roma --- ITALY  \newline
           e-mail: \ gavarini\@mat.uniroma2.it,  \newline
  \phantom{e-mail:} \hskip1pt \  gavarini\@mat.uniroma3.it
\endaddress

\abstract
    Let  $ \ghat $  be an untwisted affine Kac-Moody algebra.  The quantum
group  $ \uqg $  is known to be a quasitriangular Hopf algebra (to be
precise, a braided Hopf algebra).  Here we prove that its
unrestricted specializations at odd roots of 1 are braided too: in
particular, specializing  $ q $  at  $ 1 $  we have that the function
algebra   $ F \big[ \widehat{H} \big] $  of the Poisson proalgebraic
group   $ \widehat{H} $  dual of  $ \widehat{G} $  (a Kac-Moody group
with Lie algebra  $ \ghat $)  is braided.  This in turn implies also that
the action of the universal  $ R $--matrix  on the tensor products of
pairs of Verma modules can be specialized at odd roots of 1.
 \endabstract

\endtopmatter


\footnote""{ Keywords: {\sl Affine Quantum Groups, R-matrix}.
1991 {\it Mathematics Subject Classification:\/} 17B37, 81R50 }

\vskip0,491truecm

 \centerline{ \bf  Introduction }

\vskip11pt

\hfill \hbox{
        \vbox{
          \hbox{\it \hskip25pt   "Oh, quant'\`e affine alla sua genitrice\,! }
          \hbox{\it \hskip25pt   Osserva come anch'ella ha
belle trecce }
          \hbox{\it \hskip36pt   ch'ha ereditate dalla sua matrice" }
                     \vskip4pt
          \hbox{\sl \hskip60pt    N.~Barbecue, "Scholia" } } \hskip10pt }

\vskip11pt

  A Hopf algebra  $ H $  is called quasitriangular (cf.~[Dr], [C-P]) if there
exists an invertible element  $ R \in H \otimes H $  (or an element of
an appropriate completion of  $ H \otimes H $)  such that
  $$  \displaylines{
   {\hbox{\rm Ad}}(R) (\Delta (a)) = \Delta^{\hbox{\smallrm op}}(a)   \qquad
\forall \; a \in H  \cr
   (\Delta \otimes {\hbox{\rm id}}) (R) = R_{13} R_{23} \; ,  \qquad
({\hbox{\rm id}} \otimes \Delta) (R) = R_{13} R_{12}  \cr }  $$
where  $ \, \hbox{Ad}(R)(x) := R \cdot x \cdot R^{-1} \, $,
$ \; \Delta^{\hbox{\smallrm op}} \, $  is the opposite
comultiplication (i.~e.  $ \, \Delta^{\hbox{\smallrm op}}(a) =
\sigma \smallcirc \Delta(a) \, $  with  $ \, \sigma \colon
\, A^{\otimes 2} \to A^{\otimes 2} \, $,  $ \, a \otimes b \mapsto b
\otimes a  \, $),  and  $ \, R_{12}, R_{13}, R_{23} \in H^{\otimes 3} \, $
(or the appropriate completion of  $ H^{\otimes 3} $),  $ \, R_{12} = R
\otimes 1 \, $,  $ \, R_{23} = 1 \otimes R \, $,  $ \, R_{13} = (\sigma
\otimes {\hbox{\rm id}}) (R_{23}) =  ({\hbox{\rm id}} \otimes \sigma)
(R_{12}) \, $.
                                            \par
   As a corollary of this definition,  $ R $  satisfies the Yang-Baxter
equation in  $ H^{\otimes 3} $
  $$  R_{12} R_{13} R_{23} = R_{23} R_{13} R_{12}  $$
so that a braid group action can be defined on tensor products of
$ H $--modules  (whence applications to knot theory arise).  If  $ \ghat $
is an untwisted affine Kac-Moody algebra, the quantum universal enveloping
algebra  $ U_h \big( \ghat \big) $,  over  $ \C[[h]] $,  is quasitriangular
(cf.~[Dr], [C-P]).  On the other hand, this is not true   --- strictly
speaking ---  for its "polynomial version", the  $ \Cq $--algebra  $ \uqg $:
nonetheless, it is a braided algebra, in the sense of the following

\vskip12pt

\proclaim {Definition}  (cf. [Re1], Definition 2) A Hopf algebra  $ H $
is called braided if there exists an automorphism  $ {\Cal R}
$  of  $ H \otimes H $  (or of an appropriate completion of  $ H \otimes H
$)  distinct from  $ \sigma \colon a \otimes b \mapsto b \otimes a $  such
that
  $$  \displaylines{
   {\Cal R} \smallcirc \Delta = \Delta^{\hbox{\smallrm op}}  \cr
   (\Delta \otimes {\hbox{\rm id}}) \smallcirc {\Cal R} = {\Cal R}_{13}
\smallcirc {\Cal R}_{23} \smallcirc (\Delta \otimes{\hbox{\rm id}}) \; ,
\qquad  ({\hbox{\rm id}} \otimes \Delta) \smallcirc {\Cal R} = {\Cal R}_{13}
\smallcirc {\Cal R}_{12} \smallcirc ({\hbox{\rm id}} \otimes \Delta)
\cr }  $$
where  $ \, {\Cal R}_{12} := {\Cal R} \otimes {\hbox{\rm id}} $,  $ {\Cal R}_{23} = {\hbox{\rm id}} \otimes
{\Cal R} $,  $ {\Cal R}_{13} = (\sigma \otimes {\hbox{\rm id}}) \smallcirc ({\hbox{\rm id}}
\otimes {\Cal R}) \smallcirc (\sigma \otimes {\hbox{\rm id}}) \in Aut(H \otimes H \otimes H)\, $.
\endproclaim

\vskip12pt

   It follows from this definition that  $ {\Cal R} $  satisfies the
Yang-Baxter equation in  $ End(H^{\otimes 3}) $:
  $$  {\Cal R}_{12} \smallcirc {\Cal R}_{13} \smallcirc {\Cal R}_{23} =
{\Cal R}_{23} \smallcirc {\Cal R}_{13} \smallcirc {\Cal R}_{12}  $$
which yields a braid group action on tensor powers of  $ H $,  which is still
important for applications.  Notice that if  $ (H,R) $  is quasitriangular,
then  $ \big( H, {\hbox{\rm Ad}}(R) \big) $  is braided.
                                                    \par
   In this paper we prove that the unrestricted specializations of  $ \uqg $
at odd roots of 1 are braided too: indeed, we show that the braiding
automorphism of  $ \uqg $   --- which is, roughly speaking, the
conjugation by its universal  $ R $--matrix ---   does leave stable the
integer form   --- of   $ \uqg $  ---   which is to be "specialized".
This extends to the present case a result due to Reshetikhin (cf.~[Re1])
for the case of the quantum group  $ \, U_q\big({\frak s}{\frak
l}(2)\big) \, $,  and to Reshetikhin (cf.~[Re2]) and the author
(cf.~[Ga1]) for  $ \, U_q(\gerg) \, $,  with  $ \gerg $  finite
dimensional semisimple.  The most general case is developed in [G-H].  As
a consequence, we get that the action of the universal
$ R $--matrix
of  $ \uqg $  on tensor products of pairs of Verma modules does
specialize at odd roots of 1 as well.

\vskip17pt

\centerline { ACKNOWLEDGEMENTS }

\vskip3pt

  The author wishes to thank I.~Damiani and M.~Rosso
for many helpful discussions.

\vskip33pt

   \centerline{ \bf  \S \; 1 \,  Definitions }

\vskip15pt

   {\bf 1.1  The classical data.}  \  Let  $ \gerg $  be a simple finite
dimensional Lie algebra over the field  $ \C $  of complex numbers, and
consider the folllowing data.
                                                    \par
   The set  $ I_0 = \{1, \dots, n \} \, $,  of the vertices of the Dynkin
diagram of  $ \gerg $  (see [Bo], [Ka] for the identification between
$ I_0 $  and  $ \{1, \dots, n \} $);  a Cartan subalgebra  $ \gerh $  of
$ \gerg $;  the root system  $ \, \Phi_0 \big( \subseteq \gerh^* \big) \, $
of  $ \gerg $;  the set of simple roots  $ \, \{\, \alpha_i \mid i
\in I_0 \,\} \, $;  the Killing form  $ \, (\,\cdot\, \vert \,\cdot\, ) \, $
of  $ \gerg $,  normalized so that short roots have square length 2.  For
all  $ \, i \in I_0 \, $,  we set  $ \, d_i := {\, (\alpha_i \vert \alpha_i)
\, \over \, 2 \,} \, $.
                                                    \par
  We denote  $ \ghat $  the untwisted affine Kac-Moody algebra associated
to  $ \gerg $,  which can be realized as
$ \; \ghat = \gerg \otimes_\C \C \left[ t, t^{-1} \right] \oplus \C \cdot
c \oplus \C \cdot \partial \; $,  with the Lie bracket given by:
$ \; [c,z] = 0 \, $,  $ \, \left[ \partial, x \otimes t^m \right] =
m \, x \otimes t^m \, $,  $ \, \left[ x \otimes t^r, y
\otimes t^s \right] = [x,y] \otimes t^{r+s} + \delta_{r,-s} r \,
(x \vert y) \, c \; $  (for all  $ \, z \in \ghat $,  $ x, y \in \gerg $,
$ m \in \Z $).
                                                     \par
  For  $ \ghat $  we consider:  $ \, I := I_0 \cup \{0\}
\, $,  $ \, I_\infty := I \, \cup \{\infty\} \, $,  and
$ \, d_0 := 1 \, $;  the (generalized) Cartan matrix
$ \, A = {(a_{ij})}_{i,j \in I} \, $  (after [Ka]);  the
maximal abelian subalgebra  $ \, \hhat := \gerh \oplus \C
\cdot c \oplus \C \cdot \partial \, \big( \subseteq
\ghat \, \big) \, $;  the root system  $ \, \Phi = \Phi_+
\cup (-\Phi_+) \, \left(\, \subset {(\gerh \oplus \C \cdot c)}^* \subset
\hhat^* \, \right) \, $,  $ \, \Phi_+ = \phipre \cup
\phipim \, $  being the set of positive roots, with
$ \, \phipim = \, \{\, m \delta \mid m \in \N_+ \, \} \, $
the set of imaginary positive roots and  $ \, \phipre \, $  the
set of real positive roots.  Then  $ \ghat $  splits as
$ \, \ghat =  \hhat \oplus \big( \oplus_{\alpha \in \Phi}
\ghat_\alpha \big) \, $,  and  $ \, \text{dim}_{\C} (\ghat_\alpha) = 1 \;\,
\forall\; \alpha \in \pm \phipre \, $,  $ \; \text{dim}_{\C} (\ghat_\alpha)
= \# (I_0) = n \;\, \forall\; \alpha \in \pm \phipim \, $;
so we define the set  $ \phitildep $  of "positive roots with  multiplicity"
as  $ \, \phitildep := \phipre \cup \phitildepim \, $,  where  $ \,
\phitildepim := \phipim \times I_0 \, $.  We let  $ \, \alpha_\infty \, $
be the unique element of  $ \hhat^* $  such that  $ \, \langle
\alpha_\infty, c \rangle = 1 $,  $ \langle \alpha_\infty, \gerh \rangle
= 0 $,  $ \langle \alpha_\infty, \partial \rangle = 0 $.  Furthermore, we
set:  $ \, Q := \sum_{i \in I} \Z \cdot \alpha_i \, $,  $ \, Q_\infty :=
\sum_{i \in I_\infty} \Z \cdot \alpha_i \subset \hhat^* \, $;  for any
$ \, \beta = \sum_{i \in I} z_i \alpha_i \in Q \, $  ($ z_i \in \Z $  for
all  $ i \, $)  we set  $ \, \vert \beta \vert := \sum_{i \in I} z_i \, $.
Finally, we define the non-degenerate symmetric bilinear form on  $ \, {\Bbb R}
\otimes_\Z Q_\infty \, $  given by:  $ \, (\alpha_i \vert \alpha_j) :=
d_i a_{ij} \, (\, \forall\, i, j \in I \,) $,
$ (\alpha_\infty \vert
\alpha_j) := \delta_{0,j} \; (\, \forall\, j \in I_\infty \,) $.

\vskip12pt

{\bf 1.2  Some  $ q $-tools.} \;  For all  $ m $,  $ n $,
$ k $,  $ s \in \N_+ $,  $ n \leq m $,  we define:
$ \; {(s)}_q := {\,q^s - 1\, \over
\,q - 1\,} \, $,  $ \, {[s]}_q := {\,q^s - q^{-s}\, \over
\,q - q^{-1} \,} \, $,  $ \, {(k)}_q ! := \prod_{s=1}^k  {(s)}_q \, $,
$ \, {[k]}_q ! := \prod_{s=1}^k {[s]}_q \, $,
$ \, {\left({m \atop n} \right)}_q := {\,{(m)}_q !\, \over
\,{(m-n)}_q ! \, {(n)}_q !\,} \, $,  $ \, {\left[{m \atop n}
\right]}_q := {{[m]}_q ! \over {[m-n]}_q ! \, {[n]}_q !}
\, $  (all belonging to  $ \, \Z \big[ q, q^{-1}
\big] \, $).  For later use, we define also:
$ \; q_\alpha :=
q^{\,(\alpha \vert \alpha)\, \over \,2\,} \, $  for all
$ \, \alpha \in
\phipre \, $,  $ \; q_\alpha := q^{d_i} \, $  for all
$ \, \alpha = (r \delta, i) \in \phitildepim \, $,  $ \; q_i := q_{\alpha_i} = q^{d_i} \, $  for all  $ \, i \in I \, $.
                                             \par
  Second, we define the symbol
$ \; {(a;q)}_n := \prod_{k=0}^{n-1} (1 - a q^k) \, $,  for
$ n \in \N $,  $ a \in \C $.  Now consider the function of  $ z \, $:
$ \; {(z \, ; q)}_\infty := \prod_{n=0}^\infty (1- z q^n) \, $  \;  to be
thought of as an element of  $ \Cq[[z]] $:  if  $ q $  is a complex
number such that  $ \, |q| < 1 \, $,  the infinite product
expressing  $ {(z \, ; q)}_\infty $  converges to an  analytic function
of  $ z $  in any finite part of  $ \C $;  its  Taylor series is then
$ \; {(z \, ; q)}_\infty = \sum_{n=0}^\infty {{(-1)}^n q^{n \choose 2} \over
{(q;q)}_n} z^n \, $.  Define also  $ \; \exp_q(z) := \sum_{n=0}^\infty
{\,1\, \over \, (n)_{q^2} ! \,} \, z^n \, $;  \; then one has
  $$  \exp_q(z) = e_{q^2} \Big( \big( 1 - q^2 \big) z \Big)
= {{\Big( \big( 1 - q^2 \big) z \, ; \, q^2 \Big)}_\infty}^{\!\!\! -1} \, .  $$
   \indent   The following lemma describes the behavior of
$ \, {(z \, ; q)}_\infty \, $  for  $ \, q \rightarrow
\varepsilon \, $,   $ \, \varepsilon $  a root of 1.

\vskip12pt

\proclaim {Lemma 1.3}  ([Re1], Lemma 3.4.1; [Ga], Lemma 2.2)  Let  $ \varepsilon $  be a primitive  $ \ell $-th  root of
$ 1 $,  with  $ \ell $  odd.  The asymptotic behavior of the function (of  $ q $)  $ {(z \, ; q)}_\infty $  for  $ \, q \rightarrow \varepsilon \, $  is given by
  $$  {(z \, ; q)}_\infty = \exp \left({1 \over {q^{\ell^2}
- 1}} \sum_{n=1}^\infty {1 \over n^2} \cdot z^{\ell n} \right) \cdot {(1 - z^\ell)}^{-1/2} \cdot
\prod_{k=0}^{\ell - 1} {\left( 1 - \varepsilon^k z \right)}^{k/\ell} \cdot \big( 1 + {\Cal O}(q-\varepsilon) \big) \; .  \quad  \square  $$
\endproclaim

%
%
 \eject
   {\bf 1.4  The quantum group  $ \uqg $.}  \  The quantized  universal
enveloping algebra  $ \uqg $  (cf.~e.g.~[Dr]) is the unital associative
$ \Cq $--algebra  with generators  $ \; F_i $,  $ K_\mu $,  $ E_i \; $  ($ \, i \in I
\, $,  $ \, \mu \in Q_\infty \, $)  and relations (for all  $ \, \mu, \nu
\in Q_\infty \, $,  $ \, i, j, h \in I \, $,  $ \, i \not= j \, $)
  $$  \displaylines {
   K_\mu K_\nu = K_{\mu + \nu} = K_\nu K_\mu  \quad , \qquad K_0 = 1  \cr
   K_{\mu} E_i = q^{(\mu \vert \alpha_i)} E_i K_{\mu} \, ,
\quad  K_{\mu} F_i = q^{-(\mu \vert \alpha_i)} F_i K_{\mu} \, ,  \quad
E_i F_h - F_h E_i = \delta_{ih} {\; K_{\alpha_i} - K_{-\alpha_i} \;
\over \; q_i - q_i^{-1} \;}  \cr
   \sum_{k = 0}^{1-a_{ij}} (-1)^k {\left[ { 1-a_{ij} \atop k }
\right]}_{q_i} \!\! E_i^{1-a_{ij}-k} \! E_j E_i^k = 0 \, ,  \quad \,
\sum_{k = 0}^{1-a_{ij}} (-1)^k {\left[{1-a_{ij} \atop k} \right]}_{q_i}
\!\! F_i^{1-a_{ij}-k} \! F_j F_i^k = 0  \cr }  $$
   \indent   A Hopf algebra structure on  $ \, \uqg \, $  is defined by
($ i \in I ; \, \mu \in Q_\infty $)
  $$  \matrix
   \Delta(F_i) := F_i \otimes K_{-\alpha_i} + 1 \otimes F_i \,,  &
S(F_i):= -F_i K_{\alpha_i} \, ,  &  \epsilon(F_i) := 0 \, \phantom{.}  \\
   \Delta(K_\mu) := K_\mu \otimes K_\mu \, ,  &   S(K_\mu):= K_{-\mu} \, ,
&  \epsilon(K_\mu) := 1 \, \phantom{.}  \\
  \Delta(E_i) := E_i \otimes 1 + K_{\alpha_i} \otimes E_i \, ,  &
S(E_i):= -K_{-\alpha_i} E_i \, ,  &  \epsilon(E_i) := 0 \, .  \\
      \endmatrix  $$
   \indent   Moreover,  $ \uqg $  has a natural Hopf algebra  $ Q $--grading,
$ \, \uqg =  \oplus_{\eta \in Q} {\uqg}_\eta \, $.
                                                \par
  Let  $ \, \uqp $,  $ \uqz $,  $ \uqm $  be the subalgebras of
$ \uqg $  respectively generated by  $ \, \{\, E_i \mid i \in I \,\} $,
$ \{\, K_\nu  \mid \nu \in Q_\infty \,\} $,  $ \{\, F_i \mid i \in I
\,\} $;  then  $ \uqp $  and  $ \uqm $  are both naturally graded by  $ \, Q_+ := \sum_{i \in I} \N \cdot \alpha_i (\subset Q) \, $.  Finally, let  $ \, \uqbp := \uqp \cdot \uqz = \uqz \cdot \uqp \, $,
$ \, \uqbm := \uqm  \cdot \uqz = \uqz \cdot \uqm \, $,  to be called
{\it quantum Borel  (sub)algebras}:  these are  {\sl Hopf subalgebras}
of  $ \uqg $.

\vskip4pt

   $ \underline{\hbox{\it Remark}} \, $: \, In the definition of  $ \uqg $
several choices for the "toral part"  $ \uqz $  are possible, mainly
depending on the choice of any lattice  $ M $  such that  $ \, Q_0 \leq M
\leq P_0 \, $,  $ P_0 $  being the weight lattice of  $ \gerg $  (cf.~for
instance [B-K]).  All what follows holds as well for every such choice, up
to suitably adapting the statements involving the toral part.

\vskip12pt

   {\bf 1.5  Quantum root vectors.}  \  It is known (cf.~[Be1], [Be2])
that one can define a total ordering on  $ \phitildep $,  and accordingly
define quantum root vectors: from now on, we assume a total ordering be
fixed and quantum root vectors be defined as in [Ga2], \S~2, so that  $
\, E_\alpha \, $,  resp.~$ \, F_\alpha \, $,  is the quantum root vector
in  $ \uqp $,  resp.~in  $ \, \uqm \, $,  attached to the positive,
resp.~negative, root (with multiplicity)  $ \alpha $,  resp.~$ -\alpha $
(for any  $ \, \alpha \in \phitildep \, $).

\vskip12pt

  {\bf 1.6  Integer forms.} \  The main interest of quantum groups is to
specialize them at roots of 1: thus we need suitable integer forms of them.
                                             \par
  First, let  $ {\frak R} $  be the set of all roots of 1 (in  $ \C $)  whose
order is either 1 or an odd number  $ \ell $  with  $ \, g.c.d.(\ell,n+1)
= 1 \, $  if  $ \gerg $  is of type  $ A_n \, $,  $ \, \ell \notin 3 \N_+
\, $  if  $ \gerg $  is of type  $ E_6 $  or  $ G_2 \, $;  then let
$ \, {\Cal A} \, $  be the subset of  $ \, \C(q) \, $  of rational
functions of  $ q $  with no poles in
$ {\frak R} $.  Second, define renormalized root vectors by
$ \; \echeck_\alpha := \left( q_\alpha - q_\alpha^{-1} \right)
E_\alpha \, $,  $ \; \fcheck_\alpha := \left( q_\alpha^{-1}
- q_\alpha \right) F_\alpha \, $,  \; for all  $ \, \alpha
\in \phitildep \, $,  and let  $ \, \calUqg \, $  be the
$ {\Cal A} $--subalgebra of  $ \uqg $  generated by
$ \; \left\{ \, \fcheck_\alpha, K_\mu, \echeck_\alpha \,\big\vert\, \alpha
\in \phitildep \, , \mu \in Q_\infty \,\right\} \, $:  this is a
$ Q $--graded  Hopf subalgebra (and an  $ {\Cal A} $--form)  of
$ \uqg $  (cf.~[B-K]).  We define also
$ \; \calUqbm := \calUqg \cap \uqbm \, $,
$ \, \calUqbp := \calUqg \cap \uqbp \, $.

%
%
 \eject

\centerline {\bf \S\; 2 \  Braiding of quantum enveloping algebras }

\vskip15pt

{\bf 2.1  More notation.} \;  As we said, it is well known (cf.~[Dr]) that
the quantum algebra  $ U_h \big( \ghat \big) $  (defined over the ring
$ \C[[h]] $)  is quasitriangular; this is proved by means of Drinfeld's
method of the "quantum double".  On the other hand, for the  $ \Cq $--algebras
$ \uqg $  the correct statement is that they are braided.  To see this, we
define a suitable completion of  $ {\uqg}^{\otimes 2} $,  namely  $ \;
{\uqg}^{\widehat \otimes 2} := \left\{\, \sum_{n=0}^{+\infty} {\Cal E}_n
\cdot P_n^- \otimes P_n^+ \cdot {\Cal F}_n \,\right\} \; $  where  $ \,
P_n^- \in \uqbm $,  $ P_n^+ \in \uqbp $,  $ {\Cal E}_n \in \sum_{\vert \beta
\vert = n} {\big( \uqg \big)}_\beta $,  $ {\Cal F}_n \in \sum_{\vert \beta
\vert = -n} {\big( \uqg \big)}_\beta \, $.  It is clear that
$ {\uqg}^{\widehat \otimes 2} $  is a completion of
$ {\uqg}^{\otimes 2} $  as a Hopf algebra.
                                                  \par
  Define new quantum root vectors  $ \, \fdot_\alpha \, $  and
$ \, \edot_\alpha \, $  ($ \alpha \in \phitildep $)  as follows (like in the
proof of Proposition 4.6 in [Ga2]).  For all  $ \, \alpha \in \phipre \, $,
set  $ \, \fdot_\alpha := \fdot_\alpha \, $.  For all  $ \, \alpha =
(r \delta,i) \in \phitildepim \, $,  consider the matrix
  $$  M_r := {\left( {\big( o(i) o(j) \big)}^r \, {[a_{ij}]}_{q_i^r}
\right)}_{i,j \in I_0}  $$
where  $ \, o(i) = \pm 1 (i \in I_0) \, $  is defined in such a way
that  $ \, o(h) o(k) = -1 \, $  whenever  $ \, a_{h{}k} < 0 \, $:  then
$ \, det \big( M_r \big) \, $  is an invertible element of  $ {\Cal A} $
(see [Ga2] for the exact value), so the inverse matrix  $ \, M_r^{-1} =
{\big( \mu_{ij} \big)}_{i,j \in I_0} \, $  has all its entries in
$ {\Cal A} \, $;  now define
  $$  \fdot_{(r \delta,i)} := \sum_{j \in I_0} \mu_{ji}
F_{(r \delta,j)} \; .  $$
Similarly we define positive root vectors  $ \, \edot_\alpha \, $,
for all  $ \alpha \in \phitildep \, $.

  Now set  $ \, \exp_\alpha := \exp_{q_\alpha} \, $,  for
$ \, \alpha \in \phipre \, $,  and  $ \, \exp_\alpha :=
\exp \, $,  for  $ \, \alpha \in \phitildepim \, $;  set also
$ \, a_\alpha := 1 \, $  for  $ \, \alpha \in \phipre \, $  and
$ \, a_\alpha := {\, r \, \over \, {[r]}_q {[d_i]}_q \,} \, $  for
$ \, \alpha = (r \delta,i) \in \phitildepim \, $.

\vskip12pt

\proclaim {Theorem 2.2}  Let  $ \rcalzero $  be the algebra
automorphism of  $ {\uqg}^{\widehat \otimes 2} $  defined by
  $$  \aligned
     \rcalzero (K_\mu \otimes 1) := K_\mu \otimes 1 \; ,  &
\quad \rcalzero (1 \otimes K_\mu) := 1 \otimes K_\mu  \\
     \rcalzero(E_i \otimes 1) := E_i \otimes K_{-\alpha_i} \; ,  &  \quad
\rcalzero (1 \otimes E_i) := K_{-\alpha_i} \otimes E_i  \\
     \rcalzero (F_i \otimes 1) := F_i \otimes K_{\alpha_i} \; ,  &  \quad
\rcalzero (1 \otimes F_i) := K_{\alpha_i} \otimes F_i  \\
      \endaligned  $$
($ i \in I $,  $ \mu \in Q_\infty $)  and let  $ \, \runo
\in {\uqg}^{\widehat \otimes 2} \, $  be defined as the ordered product
  $$  \runo := \prod_{\alpha \in \phitildep} \exp_\alpha \!
\Big( a_\alpha \big( {q_\alpha}^{-1} - q_\alpha \big)
E_\alpha \otimes \fdot_\alpha \Big) = \prod_{\alpha \in
\phitildep} \exp_\alpha \! \Big( a_\alpha \big(
{q_\alpha}^{-1} - q_\alpha \big) \edot_\alpha \otimes
F_\alpha \Big)  $$
   \indent  Then  $ \, \Big( \uqg, {\hbox{\rm Ad}} \big(
\runo \big) \smallcirc \rcalzero \Big) \, $  is a braided
Hopf algebra (with  $ \runo $  as  $ R $-matrix,  in the sense of [Re1],
Definition 3).
\endproclaim

 \eject

\demo{Proof}  This is essentially proved in [Da2]: to get exactly the
present claim, one just has to take into account the following.  The formula
for the  $ R $--matrix  given in [Da2] is obtained by computing bases (of PBW
type), in quantum Borel subalgebras of opposite sign, which are orthogonal to
each other with respect to a certain perfect Hopf pairing: Lemma 2.4 in [Ga2]
extends the result of [Da2] to other possible bases, specifying which
choices of quantum root vectors are "admissible", i.e.~are such that starting
from them the construction in [Da2] still works and gives similar orthogonal
bases; finally, the remarks in the proof of Proposition 4.6 in [Ga2] show
that both the choice of quantum root vectors  $ E_\beta $  and
$ \fdot_\gamma $  ($ \beta, \gamma \in \phitildep $)  and the choice
of  $ \edot_\beta $  and  $ F_\gamma $  ($ \beta, \gamma \in \phitildep
$) are admissible (in the previous sense).   $ \square $
\enddemo

\vskip12pt

   {\bf 2.3  The braiding structure at roots of 1.} \  Our goal now is to
show that  $ \calUqg $  is braided: to be precise, we  could say that the
braiding structure of  $ \uqg $  gives by restriction a braiding structure
for  $ \calUqg $.  To begin with, we define a suitable completion of
$ \calUqg^{\otimes 2} $  (mimicking \S 2.1),  namely
  $$  \calUqg^{\widehat \otimes 2} := \left\{\,
\sum_{n=0}^{+\infty} {\overline{\Cal E}}_n \cdot P_n^-
\otimes P_n^+ \cdot {\overline{\Cal F}}_n \,\right\}  $$
where  $ \, P_n^- \in \calUqbm $,  $ P_n^+ \in \calUqbp $,
$ {\overline{\Cal E}}_n \in \sum_{\vert \beta \vert = n} {\big( \calUqg
\big)}_\beta $,  $ {\overline{\Cal F}}_n \in \sum_{\vert \beta \vert = -n}
{\big( \calUqg \big)}_\beta \, $.  It is clear that
$ \calUqg^{\widehat
\otimes 2} $  is a completion of  $ \calUqg^{\otimes 2} $  as Hopf algebra,
and that  $ \, \calUqg^{\widehat \otimes 2} \subseteq \uqg^{\widehat
\otimes 2} \, $  via the natural embedding  $ \, \calUqg \,
{\lhook\joinrel\relbar\joinrel\rightarrow} \, \uqg \, $.
                                                  \par
  Moreover, for all  $ \, \alpha \in \phitildep \, $  we
define  $ \, \fdotcheck_\alpha := \big( q_\alpha - q_\alpha^{-1} \big) \fdot_\alpha \, $,
$ \, \edotcheck_\alpha := \big( q_\alpha - q_\alpha^{-1}
\big) \edot_\alpha \in \calUqg \, $.

\vskip12pt

   For any  $ \, \varepsilon \in {\frak R} $,  we call  $ \,
{\Cal U}_\varepsilon (\ghat) \, $  the specialization
of  $ \calUqg $  at  $ \, q = \varepsilon \, $,  that is
  $$  {\Cal U}_\varepsilon (\ghat) := \, \calUqg \Big/ (q-
\varepsilon) \, \calUqg \; .  $$

\vskip12pt

\proclaim {Theorem 2.4}  The restriction of  $ \rcalzero $
(cf.~Theorem 2.2) to  $ \calUqg^{\widehat \otimes 2} $  is given by
  $$  \aligned
     \check{\Cal R}^{(0)} (K_\mu \otimes 1) = K_\mu \otimes 1 \; ,  &
\quad  \check{\Cal R}^{(0)} (1 \otimes K_\mu) = 1 \otimes K_\mu   \\
     \check{\Cal R}^{(0)} (\echeck_\alpha \otimes 1) = \echeck_\alpha \otimes
K_{-\alpha} \; ,  &  \quad   \check{\Cal R}^{(0)} (1 \otimes \echeck_\alpha) =
K_{-\alpha} \otimes \echeck_\alpha   \\
     \check{\Cal R}^{(0)} (\fcheck_\alpha \otimes 1) = \fcheck_\alpha \otimes
K_\alpha \; ,  &  \quad   \check{\Cal R}^{(0)} (1 \otimes \fcheck_\alpha) =
K_\alpha \otimes \fcheck_\alpha   \\
       \endaligned  $$
($ \, \mu \in Q_\infty, \, \alpha \in \phitildep \, $)
thus  $ \rcalzero $  restricts to an algebra automorphism
$ \check{\Cal R}^{(0)} $  of  $ \calUqg^{\widehat
\otimes 2} $.  Moreover,  let  $ \, \runo \in
\uqg^{\widehat \otimes 2} \, $  be given (as in Theorem 2.2) by
  $$  \runo := \prod_{\alpha \in \phitildep} \exp_\alpha \!
\Big( a_\alpha \big( {q_\alpha}^{-1} - q_\alpha \big)
E_\alpha \otimes \fdot_\alpha \Big) = \prod_{\alpha \in
\phitildep} \exp_\alpha \! \Big( a_\alpha \big(
{q_\alpha}^{-1} - q_\alpha \big) \edot_\alpha \otimes
F_\alpha \Big) \; .  $$
   \indent   Then the adjoint action by  $ \runo $  leaves
$ \calUqg^{\widehat \otimes 2} $  stable; thus
$ \hbox{\rm  Ad} \left( \runo \right) $  restricts to an
automorphism  $ \check{\Cal R}^{(1)} $  of
$ \calUqg^{\widehat \otimes 2} $,  and
$ \, \left( \calUqg, \check{\Cal R} \right) \, $
--- with  $ \, \check{\Cal R} := \check{\Cal R}^{(1)}
\smallcirc \check{\Cal R}^{(0)} \, $ ---   is a braided Hopf algebra.
\endproclaim

 \eject

\demo{Proof}  The first part of the statement is trivial, and the third is a
direct consequence of the first, the second, and Theorem 2.2.  To prove that
$ \hbox{\rm  Ad} \left( \runo \right) $  stabilizes  $ \calUqg^{\widehat
\otimes 2} $,  we apply an idea of Reshetikhin.
                                                \par
   We look at the different factors  $ R^{(1)}_{\,\alpha} $
(for  $ \, \alpha \in \phitildep \, $)  in the product defining  $ \runo \, $.  Notice that  $ {\uqg}^{\widehat
\otimes 2} $  has a natural  "$ Q $--pseudograding",  extending that of  $ {\uqg}^{\otimes 2} \, $:
so we can look at the homogeneous summands of  $ \runo $,  and then we find that each of them is given by the product of  {\sl finitely many}  factors  $ R^{(1)}_{\,\alpha} \, $.  Therefore, to prove the claim we have only to
show that the adjoint action by every factor
$ R^{(1)}_{\,\alpha} $  leaves
$ {\calUqg}^{\widehat \otimes 2} $  stable.
                                                \par
   For a real root  $ \, \alpha \in \phipre \, $,  the
factor  $ R^{(1)}_{\,\alpha} $  is of type
  $$  R^{(1)}_{\,\alpha} := \exp_\alpha \! \Big( a_\alpha
\big( {q_\alpha}^{-1} - q_\alpha \big) E_\alpha \otimes \fdot_\alpha \Big) = \exp_{q_\alpha} \! \Big( {\big( {q_\alpha}^{-1} - q_\alpha \big)}^{-1}
\echeck_\alpha \otimes \fcheck_\alpha \Big) \, .  $$
Using the identity in \S 1.2 this reads
  $$  R^{(1)}_{\,\alpha} = {{\Big( q_\alpha \cdot
\echeck_\alpha \otimes \fcheck_\alpha \, ; \, q_\alpha^{\,2}
\Big)}_{\!\infty}}^{\!\!\! -1} \, .  $$
Now apply Lemma 1.3: it gives
  $$  \displaylines{
   {} \qquad   R^{(1)}_{\,\alpha} = \exp\left({\,-1\, \over
\,q^{\ell^2}-1\,} \cdot {\,1\, \over \,2 d_\alpha\,} \cdot
\varphi \left( \echeck_\alpha^\ell \otimes \fcheck_\alpha^\ell \right) \right) \cdot   \hfill {}  \cr
   {} \hfill   \cdot {\left( 1 - \echeck_\alpha^\ell
\otimes \fcheck_\alpha^\ell \right)}^{-1/2} \cdot \prod_{k=0}^{\ell-1} {\big( 1 - \varepsilon^k \cdot \echeck_\alpha \otimes \fcheck_\alpha \big)}^{k/\ell} + {\Cal O}(q-\varepsilon)  \qquad {}  \cr }  $$
where we set  $ \, \varphi(z) := \sum_{n=1}^\infty {\,1\,
\over \,n^2\,} \, z^n \, $.  Thus  $ R^{(1)}_{\,\alpha} $   --- modulo a "tail" vanishing at  $ \, q = \varepsilon $  ---   contains the factor  $ \, {\left( 1 -
{\echeck_\alpha}^\ell \otimes {\fcheck_\alpha}^\ell
\right)}^{-1/2} \cdot \prod_{k=0}^{\ell-1} {\big( 1 -
\varepsilon^k \cdot {\echeck_\alpha}^\ell \otimes
{\fcheck_\alpha}^\ell \big)}^{k/\ell} \, $,  which is
"harmless" (and is trivial if  $ \, \ell = 1 $,  that is
$ \, \varepsilon = 1 \, $),  but also the factor  $ \,
\exp\left({\, 1 \, \over \, q^{\ell^2} - 1 \,} {\,1\,
\over \,2 d_\alpha\,} \, \varphi \left(
\varepsilon^{d_\alpha} {\echeck_\alpha}^\ell \otimes
{\fcheck_\alpha}^\ell \right) \! \right) $,  which has a
pole at  $ \, q = \varepsilon \, $.
                                                \par
   Here we can act as in the proof of the finite case
(cf.~[Ga1], Proposition 4.2).  Recall that  $ \, \hbox{\rm
Ad}\big( \exp(x) \big) = \exp \big( \hbox{\rm ad}(x) \big)
\, $,  where  $ \, \big( \hbox{\rm ad}(x) \big) (y) := [x,y] = x y - y x \, $.  Moreover, it is known  (see [B-K], \S 2,
or [Da1], \S 3) that the images of
$ {\echeck_\alpha}^\ell $  and of  $ {\fcheck_\alpha}^\ell $
belong to the centre of the specialized algebra  $ \, {\Cal U}_\varepsilon \big( \ghat \big) := \, \calUqg \Big/ (q -
\varepsilon) \, \calUqg \, $;  therefore  $ \, {\echeck_\alpha}^\ell \otimes {\fcheck_\alpha}^\ell \, $  belong to the centre of  $ \, {\Cal U}_\varepsilon \big( \ghat \big) \otimes {\Cal U}_\varepsilon \big( \ghat \big) \, $.  This implies that
  $$  \big[ {\echeck_\alpha}^\ell \otimes
{\fcheck_\alpha}^\ell, y \otimes z \big] \in {(q -
\varepsilon)} \cdot {\Cal U}_\varepsilon \big( \ghat \big)
\otimes {\Cal U}_\varepsilon \big( \ghat \big)  $$
hence
  $$  {(q - \varepsilon)}^{-1} \big[ {\echeck_\alpha}^\ell \otimes {\fcheck_\alpha}^\ell, y \otimes z \big] \in {\Cal
U}_\varepsilon \big( \ghat \big) \otimes {\Cal U}_\varepsilon \big( \ghat \big)  $$
and this clearly implies  $ \; \hbox{\rm Ad} \left(
R^{(1)}_{\,\alpha} \right) \left( {\calUqg}^{\widehat
\otimes 2} \right) \subseteq {\calUqg}^{\widehat \otimes 2}
\, $,  q.e.d.
                                                \par
   Now consider the factor  $ R^{(1)}_{\,\alpha} $
associated to any imaginary root  $ \, \alpha = (r \delta,i)
\in \phitildepim \, $:  by definition,
  $$  R^{(1)}_{\,\alpha} := \exp_\alpha \! \Big( a_\alpha
\big( {q_\alpha}^{-1} - q_\alpha \big) E_\alpha \otimes
\fdot_\alpha \Big) = \exp \! \left( {\, -r \, \over \,
{[r]}_q {[d_i]}_q \,} \, {\, 1 \, \over \, \big(
{q_\alpha}^{-1} - q_\alpha \big) \,} \echeck_\alpha \otimes
\fdotcheck_\alpha \right) .  $$
  In this case, we have to distinguish the cases  $ \, \ell
> 1 \, $  and  $ \, \ell = 1 \, $.
                                                   \par
  If  $ \, \ell > 1 \, $,  when  $ \, \ell \! \not\vert \, r \, $
the coefficient of  $ \, \echeck_\alpha \otimes \fdotcheck_\alpha \, $  in
the right-hand-side expression above is regular at  $ \, q = \varepsilon
\, $,  thus no problem arises.  On the other hand, if  $ \, \ell \big\vert
r \, $  that coefficient has a pole (in the factor  $ \, {{[r]}_q}^{-1}
\, $)  at  $ \, q = \varepsilon \, $;  but then (see again [B-K], \S 2,
and [Da1], \S 3) the root vectors  $ \, \echeck_\alpha \, $  and  $ \,
\fdotcheck_\alpha \, $  are again central modulo  $ (q - \varepsilon) $,
and we can conclude as in the case of real roots.
                                                   \par
  If  $ \, \ell = 1 \, $,  the coefficient of  $ \,
\echeck_\alpha \otimes \fdotcheck_\alpha \, $  has a pole at  $ \, q = 1 \, $  (in the factor  $ \, {\big( {q_\alpha}^{-1} - q_\alpha \big)}^{-1} \, $).  Now, from [B-K], \S 3, we
know that  $ \, {\Cal U}_1 (\ghat) := \, \calUqg \Big/ (q-1)
\, \calUqg \, $  is commutative, so we can apply once more the same argument than before to get that  $ \; \hbox{\rm Ad} \big( R^{(1)}_{\,\alpha} \big)
\left( {\calUqg}^{\widehat \otimes 2} \right) \subseteq
{\calUqg}^{\widehat \otimes 2} \, $.   $ \square $
\enddemo

\vskip12pt

   Let  $ \widehat{G} $  be a connected Kac-Moody group with Lie algebra
$ \ghat $,  and let  $ \widehat{H} $  be the Poisson proalgebraic group
dual of  $ \widehat{G} $  (in the sense of [B-K]): so  $ \widehat{G} $  is
a Poisson proalgebraic group whose tangent Lie bialgebra is  $ \gerg^* $.
We denote by  $ \, F \big[ \widehat{H} \, \big] \, $  the Poisson Hopf
algebra of (algebraic) regular functions on  $ \widehat{H} $.

\vskip12pt

\proclaim{Corollary 2.5}
                                            \hfill\break
   \indent   (a) \; For any  $ \, \varepsilon \in
{\frak R} $,  let  $ {\Cal R}_\varepsilon $  be the algebra
automorphism  of  $ \, {{\Cal U}_\varepsilon (\ghat)}^{\widehat \otimes 2} $  given by specialization of
$ \check{\Cal R} $  at  $ \, q = \varepsilon \, $.  Then
$ \big( {\Cal U}_\varepsilon (\ghat), {\Cal R}_\varepsilon
\big) $  is a braided Hopf algebra.
                                            \hfill\break
   \indent   (b) \; The algebra  $ \, F \big[ \widehat{H} \,
\big] \, $  is braided, by a braiding automorphism which is
one of Poisson algebra.
\endproclaim

\demo{Proof} Claim  {\it (a)}  is a direct consequence of Theorem 2.4.  As
for claim  {\it (b)},  first we recall   --- from [B-K], \S 4 ---   that
there exists a Poisson Hopf algebra isomorphism
  $$  {\Cal U}_1 (\ghat) \cong F \big[ \widehat{H} \, \big]  $$
thus the first part of claim  {\it (b)}  is nothing but a special
case of  {\it (a)}.
                                                      \par
   In addition, the Poisson bracket on  $ \, {\Cal U}_1 (\ghat) \, $
is defined, as usual, by
  $$  \{x,y\} := {\; x' y' - y' x' \; \over \; q-1 \;} {\Bigg\vert}_{q=1}
\eqno (\star)  $$
for all  $ \; x, y \in {\Cal U}_1 (\ghat) \, $,  with  $ \, x', y'
\in \calUqg \, $  such that  $ \; x = x'{\Big\vert}_{q=1} $,  $ \; y =
y'{\Big\vert}_{q=1} $;  of course a like formula defines the Poisson
bracket on the completion  $ \, {{\Cal U}_1 (\ghat)}^{\widehat \otimes 2}
\, $.  Now, since  $ \check{\Cal R} $  is an algebra automorphism
(of  $ \, {\calUqg}^{\widehat \otimes 2} \, $)  its specialization
$ \, {\Cal R}_1 \, $  automatically preserves the Poisson bracket
$ (\star) $,  i.e.~it is a Poisson algebra automorphism, q.e.d.
$ \square $
\enddemo

\vskip12pt

   $ \underline{\hbox{\it Remark}} \, $:  The results in Theorem 2.4 and
Corollary 2.5 was first proved for the finite dimensional case of the Lie
algebra  $ \, {\frak s}{\frak l}(2) \, $  in [Re1]; the case of any finite
dimensional semisimple Lie algebra was developed (and solved) in [Re2] and
in [Ga1].  The (affine) case of  $ \, \ghat = \widehat{{\frak s}{\frak l}}(2)
\, $  has been done in [H-S].  The most general situation, dealing with any
quasitriangular Lie bialgebra   --- giving rise to a quantized universal
enveloping algebra which is quasitriangular as a Hopf algebra
---   is treated in [G-H].

\vskip12pt

   {\bf 2.6  The geometry of the  $ R $--matrix  action: comparison with
the finite case.} \  In [Ga1], \S 4, the geometrical meaning of the
braiding of quantum groups of finite type at roots of 1 is explained.  The
main points are Theorems 4.4--5, where one shows that the braiding
automorphism  $ {\Cal R}_1 $  (in the present notation) is more than a formal
object   --- defined on some completion made of formal series ---   for it
maps rational functions (on the dual Poisson group times itself) onto rational
functions: hence it defines a birational automorphism of the square of this dual group (as a
complex variety), which enjoys nice properties.  The key step in the proof of
such a result exploits the fact that the Hamiltonian vector fields associated
to the functions  $ \echeck_\alpha $  and  $ \fcheck_\alpha $,  for all
positive roots  $ \alpha $  (real, since we are in the finite case), are
integrable (we consider root vectors  $ \echeck_\alpha $,  $ \fcheck_\alpha $
at  $ \, q = 1 \, $  as holomorphic functions on the dual Poisson group
$ \, H \cong Spec \big( {\Cal U}_1 (\gerg) \big) \, $).
                                             \par
   In the affine case, the situation is in part similar: one can treat the
factors of the  $ R $--matrix  associated to real positive roots exactly as
in the finite case (compare the first part of the proof of Theorem 2.4 above
with the proof of Proposition 3.7 in [Ga1]), and everything works (as in the
proof of Proposition 4.2 in [Ga1]) because the Hamiltonian vector fields which
occur   --- associated to real root vectors ---   are again integrable; but
in the case of imaginary positive roots one has to deal with Hamiltonian
vector fields   --- associated to  $ \echeck_\alpha $  and
$ \fdotcheck_\alpha $,  for imaginary  $ \alpha $  ---   which  {\sl are not
integrable}.
                                             \par
   So in the affine case the most one gets is that the
braiding defines an automorphism (with ``nice'' properties)
of the dual  {\it formal}  Poisson group associated to  $ \,
\widehat{H} \times \widehat{H} \, $.  This last result
extends, via a different approach, to the general case of
any quasitriangular Lie bialgebra: see [G-H].

\vskip12pt

   {\bf 2.7  The  $ R $--matrix  action on Verma modules.} \  For any commutative unital ring  $ A $,  denote by
$ A^\star $  the group of invertible elements of  $ A $.
Given  $ \, \lambda = {(\lambda_i)}_{i \in I_\infty} \in
{\big( {\Cq}^\star \big)}^{n+2} \, $,  let  $ V_q(\lambda) $
be the Verma module (for  $ \uqg $) of highest weight
$ \lambda $.  We recall that it is defined as follows:
define on the line  $ \, \Cq \cdot v_\lambda \, $  a
structure of  $ \uqbp $--module  by
  $$  E_i.v_\lambda := 0 \, ,  \quad  K_j.v_\lambda := \lambda_j v_\lambda   \qquad \quad  \forall \, i \in I,  j \in I_\infty \; ;  $$
then  $ V_q(\lambda) $  is by definition the  $ \uqg $--module  induced by  $ \Cq.v_\lambda $;  in particular, it is
a free  $ \uqm $ --module of rank 1, hence it is  $ Q_+ $--graded:  $ \, V_q(\lambda) = \oplus_{\eta \in Q_+} {\big( V_q(\lambda) \big)}_\eta \, $,  with  $ \, K_i.v = \lambda_i q^{-(\alpha_i \vert \eta)} \cdot v \, $  for all  $ \, i \in I_\infty $,  $ \, v \in {\big( V_q(\lambda) \big)}_\eta
\, $,  $ \eta \in Q_+ \, $.
                                                  \par
   Now assume  $ \, \lambda = {(\lambda_i)}_{i \in I_\infty} \in {\big( {\Cal A}^\star \big)}^{n+2} \, $:  then
$ V_q(\lambda) $  is also an  $ \calUqg $--module,  and \
$ \, {\Cal V}_q(\lambda) := \calUqg.v_\lambda \, $  is an
$ \calUqg $--module.  It is also clear that  $ \, {\Cal V}_q(\lambda) \, $  is a free  $ \calUqm $ --module of rank 1, and it is of course  $ Q_+ $--graded  as well.
                                                  \par
   For any  $ \, \varepsilon \in {\frak R} \, $,  we denote  $ {\Cal V}_\varepsilon (\lambda) $  the specialization of
$ {\Cal V}_q (\lambda) $  at  $ \, q = \varepsilon $,  i.e.
  $$  {\Cal V}_\varepsilon (\lambda) := {\Cal V}_q(\lambda) \Big/ (q-\varepsilon) {\Cal V}_q(\lambda) \, .  $$
   \indent   Consider on the Cartan subalgebra  $ \hhat $  the Killing form   --- which is dual of the form  $ (\ \vert \ ) $  on  $ \hhat^* $  defined in \S 1.1 ---   and let
$ T $  be its canonical element: i.e.,  $ \, T = \sum_{i \in I_\infty} u_i \otimes w_i \, $  where  $ {\{u_i\}}_{i \in I_\infty} $  and  $ {\{w_i\}}_{i \in I_\infty} $  are basis of  $ \hhat $  dual of each other with respect to the Killing form.  Let  $ \, \lambda $,  $ \mu \in {\big( {\Cal A}^\star \big)}^{n+2} \, $:  for simplicity, we assume
$ \lambda $  and  $ \mu $  to be of the form
  $$  \lambda = {\big( q^{l_i} \big)}_{i \in I_\infty} \, ,  \qquad \qquad  \;  \mu = {\big( q^{m_i} \big)}_{i \in
I_\infty}  $$
for some integers  $ l_i $,  $ m_i $  ($ i \in I_\infty $),  and we set  $ \, l := \sum_{i \in I_\infty} l_i \omega_i
\, $,  $ \, m := \sum_{i \in I_\infty} m_i \omega_i \, $,
where  $ \, {\{\omega_i\}}_{i \in I_\infty} \, $  is a basis of  $ \hhat^\ast $  dual of  $ \, {\{\alpha_i\}}_{i \in
I_\infty} \, $  with respect to the Killing form (i.e.~$ \,
(\alpha_i \vert \omega_j) = \delta_{i{}j} \;\, \forall \, i,
j \in I_\infty \, $).  We define a linear operator
  $$  \displaylines{
   q^{-T} \colon \, {\Cal V}_q (\lambda) \otimes {\Cal V}_q
(\mu) \longrightarrow {\Cal V}_q (\lambda) \otimes {\Cal
V}_q (\mu)  \cr
   \hbox{by}  \quad \qquad  q^{-T}. \big(
v'\otimes v'' \big) := q^{-(l-\eta \vert m-\xi)} \cdot
v'\otimes v'' \, ,  \qquad  \forall \, v' \in {\big( {\Cal V}_q (\lambda) \big)}_\eta \, , \;  v'' \in {\big( {\Cal
V}_q (\lambda) \big)}_\xi \, .  \cr }  $$
   \indent   For any pair of Verma modules  $ V_q(\lambda) $  and  $ V_q (\mu) $,  the algebra  $ \, \uqg^{\widehat \otimes 2} \, $  acts on  $ \, V_q(\lambda) \otimes V_q(\mu) \, $:  in fact since  $ V_q(\lambda) $  is highest weight,
it is clear that only finitely many summands in the
expansion of any element of  $ \, \uqg^{\widehat \otimes 2} \, $  act non-trivially; similarly, the algebra  $ \, \calUqg^{\widehat \otimes 2} \, $  acts on  $ \, {\Cal V}_q(\lambda) \otimes {\Cal V}_q(\mu) \, $.  As a consequence,  $ \runo $  acts as a well-defined operator on  $ \, V_q(\lambda) \otimes V_q(\mu) \, $
                                               \par
  We call  {\it universal  $ R $--matrix}  of
$ \uqg $  the formal element
  $$  \eqalign{
   R:= \runo \cdot q^{-T}  &  = \prod_{\alpha \in
\phitildep} \exp_\alpha \! \Big( a_\alpha \big(
{q_\alpha}^{-1} - q_\alpha \big) E_\alpha \otimes
\fdot_\alpha \Big) \cdot q^{-T} =  \cr
                       &  = \prod_{\alpha \in
\phitildep} \exp_\alpha \! \Big( a_\alpha \big(
{q_\alpha}^{-1} - q_\alpha \big) \edot_\alpha \otimes
F_\alpha \Big) \cdot q^{-T} \, ;  \cr }  $$
this is a universal  $ R $--matrix  for  $ \uqg $  in the sense of [Da2].

\vskip4pt

   $ \underline{\hbox{\it Remark}} \, $: \, notice that, if we deal with
the quantum group  $ U_h (\ghat) $  over the ring  $ \C[[h]] $,  the
$ R $--matrix  takes the simpler form
  $$  R := \prod_{\alpha \in \phitildep} \exp_\alpha \!
\Big( a_\alpha \big( \exp^{-d_\alpha h} - \exp^{d_\alpha}
\big) E_\alpha \otimes \fdot_\alpha \Big) \cdot
\exp (-h T)  $$
which is an element of the topological ($ h $--adically  complete)  tensor product  $ U_h (\ghat) \otimes U_h
(\ghat) $.

\vskip4pt

  For any pair of Verma modules  $ V_q(\lambda) $  and
$ V_q (\mu) $,  the  $ R $--matrix  acts as a well-defined operator on  $ \, V_q(\lambda) \otimes V_q(\mu) \, $.  Our previous results tell us that this action can be specialized at roots of 1.

\vskip12pt

\proclaim{Theorem 2.8}  The action of the universal  $ R $--matrix  on  $ \, V_q(\lambda) \otimes V_q(\mu) \, $  restricts to an action on  $ \, {\Cal V}_q (\lambda) \otimes {\Cal V}_q (\mu) \, $,  hence it specializes to an action on
$ \, {\Cal V}_\varepsilon (\lambda) \otimes {\Cal V}_\varepsilon (\mu) \, $  for any  $ \, \varepsilon
\in {\frak R} \, $.
\endproclaim

\demo{Proof}  The second part of the claim is a direct consequence of the first.
                                               \par
  Since  $ \, {\Cal V}_q(\lambda) = \calUqg.v_\lambda \, $
and  $ \, {\Cal V}_q(\mu) = \calUqg.v_\mu \, $,  we just
need to look at elements such as  $ \, R(x.v_\lambda \otimes
y.v_\mu) \, $  with  $ \, x $,  $ y \in \calUqg \, $.  We
have
  $$  \displaylines {
   {} \qquad   R(x.v_\lambda \otimes y.v_\mu) = R \big( (x \otimes y).(v_\lambda \otimes v_\mu) \big) =   \hfill  \cr
   {}   \hfill  = \big( R (x \otimes y) R^{-1} \big) \cdot
R.(v_\lambda \otimes v_\mu) = \hbox{Ad}(R) \big( x \otimes y
\big).\big(R.(v_\lambda \otimes v_\mu)\big)  \qquad {}
\cr }  $$
where  $ \, R^{-1} \, $  denotes a formal inverse to  $ R $,
(which induces the inverse operator on tensor product of
Vermas modules).  Now, definitions are given in such a way
that  $ \, \hbox{Ad}(R) \, $  coincides with the braiding
automorphism  $ \, \hbox{Ad}\left( R^{(1)} \right)
\smallcirc \rcalzero \, $  of  $ \uqg $  (cf.~Theorem 2.2,
and the  {\it Remark\/}  above): then by Theorem 2.4 we get
that  $ \, \hbox{Ad}(R) \big( x \otimes y \big) \in \calUqg
\, $,  so we only need to show that  $ \, R.(v_\lambda
\otimes v_\mu) \in {\Cal V}_q(\lambda) \otimes {\Cal
V}_q(\mu) \, $.  It is clear that  $ \, q^{-T} (v_\lambda \otimes v_\mu) = q^{(l \vert m)} \cdot v_\lambda \otimes v_\mu \in {\Cal V}_q(\lambda) \otimes {\Cal V}_q(\mu) \, $.
Moreover, from  $ \, E_i.v_\lambda = 0 \, $  for all  $ \, i
\in I \, $  we have also  $ \, E_\alpha.v_\lambda = 0 \, $
for all  $ \, \alpha \in \phitildep \, $:  then by
definition of  $ \runo $  we have  $ \, \runo.(v_\lambda
\otimes v_\mu) = v_\lambda \otimes v_\mu \, $  in  $ \,
V_q(\lambda) \otimes V_q(\mu) \, $,  hence also in  $ \,
{\Cal V}_q(\lambda) \otimes {\Cal V}_q(\mu) \, $,  q.e.d.
$ \square $
\enddemo

\vskip12pt

   {\bf Remarks 2.9:\/} \, {\it (a)} \,  As it is clear from
the proof, the previous result holds as well for  {\sl
lowest weight}  modules, and even for pair of modules in
which only the first one is highest weight or the second is
lowest weight.
                                                     \par
   {\it (b)} \,  The analogues of Theorem 2.4 and Corollary 2.5 also hold for  {\sl finite}  type quantum groups (cf.~[Ga1], \S\S 3--4): therefore Theorem 2.8 holds as well in the finite case (with the same proof).  In the case of
$ \, \gerg = {\frak s}{\frak l}(2) \, $,  such a result is complementary to another one   --- due to Date  {\it et al.}, cf.~[D-J-M-M] and [C-P], Proposition 11.1.17 ---   which concern cyclic (or periodic) representations.

\vskip1truecm

\enddocument

\bye
\end